\documentclass[journal,9pt]{IEEEtran}
\ifCLASSINFOpdf
\else
\fi
\usepackage[utf8]{inputenc}
\usepackage{graphicx}
\usepackage{amsmath} 
\usepackage{mathtools} 
\usepackage{amssymb}  
\usepackage{psfrag}
\usepackage{mathrsfs}
\usepackage{color}
\usepackage{cases}
\usepackage{pgfplots}
\usepackage{tikz}
\usetikzlibrary{arrows,automata}
\usetikzlibrary{calc}
\usetikzlibrary{intersections, positioning,fit,arrows.meta,backgrounds}
\usepackage{pgf-extrect}
\usepackage[authormarkup=none]{changes}
\usepackage[hidelinks]{hyperref}

\usepackage{microtype}

\usepackage{bm} 

\newtheorem{thm}{Theorem}
\newtheorem{Lemma}{Lemma}
\newtheorem{rem}{Remark}

\newcommand{\col}[1]{\operatorname{col}(#1)}
\newcommand{\diag}[1]{\operatorname{diag}(#1)}

\renewcommand{\d}{\mathrm{d}}
\newcommand{\e}{\mathrm{e}}

\newcommand{\ie}{i.\,e., }

\newcommand{\wrt}{w.\,r.\,t. }

\renewcommand{\t}{^{\top}}

\definechangesauthor[name=Simon,color=blue]{SK} 

\makeatletter
	\def\parsenode[#1]#2\pgf@nil{%
	    \tikzset{label node/.style={#1}}
	    \def\nodetext{#2}
	}
	
	\tikzset{
	    add node at x/.style 2 args={
	        name path global=plot line,
	        /pgfplots/execute at end plot visualization/.append={
	                \begingroup
	                \@ifnextchar[{\parsenode}{\parsenode[]}#2\pgf@nil
	            \path [name path global = position line #1-1]
	                ({axis cs:#1,0}|-{rel axis cs:0,0}) --
	                ({axis cs:#1,0}|-{rel axis cs:0,1});
	            \path [xshift=1pt, name path global = position line #1-2]
	                ({axis cs:#1,0}|-{rel axis cs:0,0}) --
	                ({axis cs:#1,0}|-{rel axis cs:0,1});
	            \path [
	                name intersections={
	                    of={plot line and position line #1-1},
	                    name=left intersection
	                },
	                name intersections={
	                    of={plot line and position line #1-2},
	                    name=right intersection
	                },
	                label node/.append style={pos=1}
	            ] (left intersection-1) -- (right intersection-1)
	            node [label node]{\nodetext};
	            \endgroup
	        }
	    },
	    add node at y/.style 2 args={
	        name path global=plot line,
	        /pgfplots/execute at end plot visualization/.append={
	                \begingroup
	                \@ifnextchar[{\parsenode}{\parsenode[]}#2\pgf@nil
	            \path [name path global = position line #1-1]
	                ({axis cs:0,#1}-|{rel axis cs:0,0}) --
	                ({axis cs:0,#1}-|{rel axis cs:1,1});
	            \path [yshift=1pt, name path global = position line #1-2]
	                ({axis cs:0,#1}-|{rel axis cs:0,0}) --
	                ({axis cs:0,#1}-|{rel axis cs:1,1});
	            \path [
	                name intersections={
	                    of={plot line and position line #1-1},
	                    name=left intersection
	                },
	                name intersections={
	                    of={plot line and position line #1-2},
	                    name=right intersection
	                },
	                label node/.append style={pos=1}
	            ] (left intersection-1) -- (right intersection-1)
	            node [label node] {\nodetext};
	            \endgroup
	        }
	    }
	}
\makeatother

\allowdisplaybreaks[1]

\begin{document}
%
\title{Robust Cooperative Output Regulation for a Network of Parabolic PDE Systems}
%
%
%

\author{Joachim~Deutscher,~\IEEEmembership{Member,~IEEE}
\thanks{J. Deutscher is with the Institut für Mess-, Regel- und Mikrotechnik, Universit\"at
    Ulm, Albert-Einstein-Allee 41, D-89081~Ulm, Germany.
    (e-mail: joachim.deutscher@uni-ulm.de)}}

%
%

\markboth{IEEE TRANSACTIONS ON AUTOMATIC CONTROL,~Vol.~XX, No.~Y, AUGUST~2020}%
{Shell \MakeLowercase{\textit{et al.}}: Bare Demo of IEEEtran.cls for Journals}
%



\maketitle

\begin{abstract}
This paper considers the robust cooperative output regulation for a network of parabolic PDE systems. The solution of this problem is obtained by extending the cooperative internal model principle from finite to infinite dimensions. For a time-invariant digraph describing the communication topology, a two-step backstepping approach is presented to systematically design cooperative state feedback regulators. They allow to solve  both the leader-follower and the leaderless output synchronization problem in the presence of disturbances and model uncertainty for a finite-dimensional leader. Solvability conditions of the robust cooperative output regulation problem are presented in terms of the communication graph and the agent transfer behaviour. The results of the paper are demonstrated for a MAS consisting of four uncertain parabolic agents with and without a finite-dimensional leader in the presence of disturbances.

\end{abstract}


\section{Introduction}\label{sec:intro}
The \emph{networked control of multi-agent systems (MAS)} has been a very active research topic for about two decades with a still increasing research effort in the control community. The interest in networked control arises from the fact that advances in communication technology by means of digital networks allow an efficient information exchange between different spatially separated systems. With this, cooperative control tasks can be systematically solved.
An overview of the state of the art as well as corresponding applications can be found in the recent monographs \cite{Bu19,Lu19,Ch19}.

Starting with simple integrators describing the agents increasingly complex system dynamics were taken into account in networked control. Currently, systematic design methods are available for linear systems, which can be found in \cite{Bu19,Lu19,Ch19}. Further advances include the extension to nonlinear systems (see, e.\:g., \cite{So19}) and to fractional systems (see, e.\:g., \cite[Ch. 7]{Re11}). In many applications it is required to take both the temporal and spatial system dynamics into account. Therefore, it is also of interest to design networked controllers for MAS with  distributed-parameter agents. Applications include industrial furnaces consisting of a network of heaters (see  \cite{Ca14}), networks of HVAC systems in building climate control (see \cite{Sa17}), networks of Lithium-Ion cells in battery management (see \cite{Qu16})  or consensus control in environmental applications (see \cite{Tr12}). Different from other system classes the networked control of distributed-parameter MAS is still an emerging research topic. Recent contributions consider parabolic agents in \cite{Dem13,Pi16} and parabolic PDEs with a diffusive coupling in \cite{Wu12,Wu16}, while networks of wave equations are investigated in \cite{Ag20,Che20}. General classes of distributed-parameter agents were dealt with in \cite{Dem18} by making use of an abstract setting. 

Common to all previous contributions for distributed-parameter MAS is the fact that the considered network is \emph{homogeneous}, i.\:e., all agents are identical. This, of course, requires some approximations, because parameter variations result from different agent environments or variations in the production process. Therefore, it is of interest to design synchronizing networked controllers, which are able to tolerate at least sufficiently small parameter perturbations. A general approach to deal with the corresponding synchronization problem is \emph{cooperative output regulation} (see, e.\:g., \cite{Hua12,Wu17,Li19}), which generalizes the classical output regulation problem (see, e.\:g., \cite{Hu04} for lumped-parameter systems). In particular, the leader-follower output synchronization can be seen as an extension of the output regulation theory, in which the reference model plays the role of the leader and several plants are the agents. Therein, not all agents have access to the reference input due to the communication constraints. Consequently, synchronization can only be achieved cooperatively by an information exchange between the agents through the communication network. Different from the usual \emph{leader-follower synchronization problem}, however, also disturbances are taken into account and the leader may differ from the followers. Furthermore, by omitting the reference model, i.\:e., the leader, also the \emph{leaderless output synchronization problem} is contained as a special case. This was soon recognized for finite-dimensional systems in the works \cite{Wi11,Su12a,Su12} and led to the development of a \emph{distributed} or \emph{cooperative observer}, which allows to solve the cooperative output regulation problem by feedforward control. For the latter the so-called regulator equations have to be solved. Since they depend on system parameters, this approach is not robust. In contrast, by assuming that the outputs to be synchronized are available for measurement, the \emph{distributed} or \emph{cooperative internal model principle} can be applied to achieve cooperative output regulation in the presence of non-destabilizing parameter perturbations (see \cite{Wa10,Su13}). For this, it is not required to solve the regulator equations. The recent monographs \cite{Li19,Wu17} demonstrate that cooperative output regulation is still a very active research area for linear finite-dimensional MAS, where also a more detailed literature overview can be found. 
Recent advances in the backstepping-based solution of the output regulation problem for distributed-parameter systems (see, e.\:g., \cite{Deu19,Deu19a} for parabolic systems and \cite{An16,Deu17a,Deu17b} for hyperbolic systems as well as the references in these works) suggest to generalize these results to networks of parabolic PDEs, in order to solve the cooperative robust output regulation problem also for distributed-parameter MAS.

This paper considers the \emph{robust cooperative output regulation problem} for MAS with boundary controlled agents subject to spatially varying coefficients. The outputs of the agents can be defined in-domain pointwise and distributed as well as at the boundaries. The reference input coinciding with the output of the leader (reference model) and the local disturbances affecting each agent in-domain, at the boundaries and at the output are generated by a finite-dimensional signal model. 
In order to provide a systematic solution for parabolic agents by making use of the cooperative internal model principle, the results from \cite{Wa10,Su13} for the robust cooperative output regulation problem in finite dimensions are combined with the results in \cite{Deu19} for the robust output regulation in infinite dimensions. This requires to add a cooperative internal model to the plant and to stabilize the resulting augmented system. For this, the communication topology has to be taken into account, which is described by a time-invariant digraph. It is shown that this stabilization problem is solvable  for a nominal homogeneous  MAS with parabolic agents if the digraph is connected and the parabolic agents satisfy the conditions for the usual output regulation (see \cite{Deu19}). Different from the latter result, however, the solution of the synchronization problem requires to solve a simultaneous stabilization problem, which is very challenging for distributed-parameter systems. Furthermore, the state feedback regulator for each agent can only use information of the agent and its neighbors resulting in constraints for the state feedback design. In order to provide a systematic design procedure for this stabilization problem, a two step backstepping approach is presented to map the closed-loop system into an exponentially stable ODE-PDE cascade with a prescribed decay rate. This is achieved by a local backstepping transformation for each agent and a cooperative decoupling transformation, which takes the communication between the agents into account. As a result, the simultaneous stabilization problem required for the output synchronization of the agents has only to be solved for the ODE subsystems. Hence, systematic methods from the literature become available for their stabilization. In the presence of model uncertainty, an uncertain heterogeneous MAS results, for which cooperative output regulation is verified provided that the nominal networked controller also stabilizes the uncertain MAS. The latter is guaranteed by the structural stabilizing property of the networked controller in the sense that the uncertain networked controlled  MAS remains stable for sufficiently small parameter perturbations. These results are subsequently used to solve the robust output synchronization problem without a leader (i.\:e., without a reference model). 

This paper demonstrates for the first time that output regulation theory provides a systematic framework to extend methods for networked controlled MAS with lumped-parameter agents to distributed-parameter agents. Furthermore, it shows that the backstepping approach (see, e.\:g., \cite{Kr08}) is also an useful tool to systematically design networked controller for distributed-parameter MAS by combining it with graph-theoretic methods.

After the problem formulation in the next section, the design of the cooperative state feedback regulator is presented in Section \ref{sec:nomout}. Robust cooperative output regulation of the resulting networked controller is verified in Section \ref{sec:robout}. These results are extended in Section \ref{sec:roboutsyncll} for the solution of the robust leaderless output synchronization problem. A MAS of four unstable parabolic agents demonstrates the results of the paper with and without an ODE-leader in the presence of disturbances and model uncertainty. 

\subsection{Elements from Graph Theory}\label{sec:gtheo}
The communication topology between the agents is described by a \emph{time invariant (weighted) digraph} $\mathcal{G}$. This is a triple $\mathcal{G} = \{\mathcal{V},\mathcal{E},A_{\mathcal{G}}\}$, in which $\mathcal{V}$ is a set of $N$ \emph{nodes} $\mathcal{V} = \{\nu_1,\ldots,\nu_N\}$, one for each agent and $\mathcal{E} \subset \mathcal{V} \times \mathcal{V}$ is a set of \emph{edges} that models the information flow from the node $\nu_j$ to $\nu_i$ with $(\nu_i,\nu_j) \in \mathcal{E}$. This flow is weighted by $a_{ij} \geq 0$, which are the element of the \emph{adjacency matrix} $A_{\mathcal{G}} = [a_{ij}] \in \mathbb{R}^{N \times N}$ with $a_{ii} = 0$, $i = 1,\ldots,N$. From this, the \emph{Laplacian matrix} $L_{\mathcal{G}} \in \mathbb{R}^{N \times N}$ of the graph $\mathcal{G}$ can be derived by $L_{\mathcal{G}} = D_{\mathcal{G}} - A_{\mathcal{G}}$, where $D_{\mathcal{G}} = \operatorname{diag}(d_1,\ldots,d_N)$ with $d_i = \sum_{k=1}^Na_{ik}$, $i = 1,\ldots,N$, is the \emph{degree matrix} of $\mathcal{G}$. Hence, the elements $l_{ij}$ of $L_{\mathcal{G}}$ are $l_{ii} =  \sum_{k=1}^{N}a_{ik}$ and $l_{ij} = -a_{ij}$, 
$i \neq j$. A \emph{path} from the node $\nu_j$ to the node $\nu_i$ is a sequence of $r \geq 2$ distinct nodes $\{\nu_{l_1},\ldots,\nu_{l_r}\}$ with $\nu_{l_1} = \nu_j$ and $\nu_{l_r} = \nu_i$ such that $(\nu_k,\nu_{k+1}) \in \mathcal{E}$. A graph $\mathcal{G}$ is said to be \emph{connected} if there is a node $\nu$, called the \emph{root}, such that, for any node $\nu_i \in \mathcal{V} \setminus \{\nu\}$, there is a path from $\nu$ to $\nu_i$. For further details on graph theory see, e.\:g., \cite[Ch. 2]{Mes10}. In the paper the \emph{Kronecker product} $A \otimes B = [a_{ij}B] \in \mathbb{C}^{n_1n_2 \times m_1m_2}$ of two matrices $A = [a_{ij}] \in \mathbb{C}^{n_1 \times m_1}$ and $B \in \mathbb{C}^{n_2 \times m_2}$ is utilized (see, e.\:g., \cite{St91}). In what follows the multiplication property $(A \otimes B)(C \otimes D) = AC \otimes BD$ with $A \in \mathbb{R}^{n_1 \times m_1}$, $B \in \mathbb{R}^{n_2 \times m_2}$, $C \in \mathbb{R}^{m_1 \times m_3}$ and $D \in \mathbb{R}^{m_2 \times m_4}$ of the Kronecker product is needed (see \cite[Ch 1.3]{St91}).


\section{Problem formulation}\label{sec:probform}
Consider a \emph{multi-agent system (MAS)} consisting of the $N > 1$ \emph{heterogeneous  parabolic agents}
\begin{subequations}\label{drs}
	\begin{align}
	\dot{x}_i(z,t) &= \bar{\lambda}_i(z)x''_i(z,t)  +  \bar{a}_i(z)x_i(z,t)  +  g_{1,i}\t(z)d_i(t)\label{pdes}\\
	x'_i(0,t) &= \bar{q}_{0,i}x_i(0,t)  +   g_{2,i}\t d_i(t), &&\hspace{-2cm}  
	t > 0\label{bc1}\\
	x'_i(1,t) &= \bar{q}_{1,i}x_i(1,t) + u_i(t) + g_{3,i}\t d_i(t),&&\hspace{-2cm}  t > 0 \label{bc2}\\
	y_i(t) &= \bar{\mathcal{C}}_i[x_i(t)] + g_{4,i}\t d_i(t),&&\hspace{-2cm}  t \geq 0\label{outcontr}
	\end{align}
\end{subequations}	
for $i = 1,\ldots,N$. The state $x_i(z,t) \in \mathbb{R}$ of \eqref{drs} is defined on $(z,t) \in (0,1) \times \mathbb{R}^+$, $\bar{\lambda}_i = 1 + \Delta\lambda_i \in C^2[0,1]$ and $\bar{a}_i = a + \Delta a_i \in C[0,1]$ are assumed. The input locations of the \emph{disturbance} $d_i(t) \in \mathbb{R}^{m_i}$ are characterized by a vector function $g_{1,i}$ with piecewise continuous elements and $g_{k,i} \in \mathbb{R}^{m_i}$, $k = 2,3,4$, which have not to be available for the controller design. In \eqref{bc1} and \eqref{bc2} the coefficients $\bar{q}_{k,i} = q_k + \Delta q_{k,i} \in \mathbb{R}$, $k = 0,1$, specify Robin or Neumann BCs, the input is $u_i(t) \in \mathbb{R}$ and the initial condition (IC) of the system reads $x_i(z,0) = x_{0,i}(z) \in \mathbb{R}$. 
\begin{rem}
Note that the assumption $\bar{\lambda}_i = 1 + \Delta\lambda_i$ and a missing advection term in \eqref{pdes} mean no loss of generality, since this can always be ensured by applying the transformations in \cite[Ch. 2.1 \& 4.8]{Kr08}. \phantom{leer} \hfill $\triangleleft$
\end{rem}

The output $y_i(t) \in \mathbb{R}$ to be controlled can be defined distributed in-domain, point-wise in-domain at $l$ locations, at the boundaries and combinations thereof. This leads to the formal \emph{output operator} 
\begin{equation}\label{cdefindomain}
\bar{\mathcal{C}}_i[h] = \int_0^1\bar{c}_i(\zeta)h(\zeta)\d \zeta + \bar{c}_{b0,i} h(0) + \bar{c}_{b1,i} h(1)
\end{equation}
for $h(z) \in \mathbb{R}$ with $\bar{c}_{bk,i} = {c}_{bk} + \Delta {c}_{bk,i}  \in \mathbb{R}, k=0,1$, and $\bar{c}_i(z) = \bar{c}_{0,i}(z) + \sum_{k=1}^l\bar{c}_{k,i}\delta(z-z_k)$, in which $\bar{c}_{0,i}(z) = c_{0}(z) + \Delta c_{0,i}(z)  \in \mathbb{R}$ with $\bar{c}_{0,i}(z)$ piecewise continuous functions, $\bar{c}_{k,i} = c_{k} + \Delta c_{k,i} \in \mathbb{R}$ and $z_{k} \in (0,1)$, $k = 1,\ldots,l$. The known \emph{nominal parameters} are $a$, $q_k$, $c_0$, $c_{bk}$, $k = 0,1$ and $c_{k}$, $z_k$, $k = 1,\ldots,l$, whereas 
\begin{equation}\label{modunc}
 \Delta \lambda_i(z), \Delta a_i(z), \Delta q_{k,i}, \Delta c_{bk,i}(z), \Delta c_{0,i}(z)   \text{ and } \Delta c_{k,i}
\end{equation}
represent unknown \emph{model uncertainties}. With this, the nominal agents give rise to a \emph{homogeneous MAS}.

For all agents a common \emph{reference input} $r(t) \in \mathbb{R}$ is specified by the solution of the \emph{global reference model}
\begin{subequations}\label{smodelg}
	\begin{align}
	\dot{w}_r(t) &= S_rw_r(t), && t > 0, \quad w_r(0) = w_{r,0} \in \mathbb{R}^{n_r} \label{sigmodssg}\\
	r(t) &= p_r\t w_r(t), &&  t \geq 0\label{sigmoy_routg}
	\end{align}
\end{subequations}
with $p_r \in \mathbb{R}^{n_r}$ and the pair $(p_r\t,S_r)$ observable. It is assumed that the \emph{spectrum} $\sigma(S_r)$ of $S_r \in \mathbb{R}^{n_r \times n_r}$ has only eigenvalues on the imaginary axis, i.\:e., $\sigma(S_r) \subset \text{j}\mathbb{R}$, and that $S_r$ is diagonalizable. Hence, \eqref{smodelg} describes a wide class of reference inputs including constant and  trigonometric functions of time as well as linear combinations thereof. The extension to non-diagonalizable global reference models is also possible by making use of the results in \cite{Deu19}. The disturbances $d_i$, $i = 1,\ldots,N$, acting on the individual agents are described by the \emph{local disturbance models}
\begin{subequations}\label{smodel}
	\begin{align}
	\dot{w}_{d_i}(t) &= S_{d_i}w_{d_i}(t), &&\hspace{-0.2cm} t > 0, \quad w_{d_i}(0) = w_{d_i,0} \in \mathbb{R}^{n_{d_i}} \label{sigmodss}\\
	d_i(t) &= P_{d_i}w_{d_i}(t), &&\hspace{-0.2cm}  t \geq 0,\label{sigmoy_rout}
	\end{align}
\end{subequations}
in which $P_{d_i} \in \mathbb{R}^{m_i \times n_{d_i}}$, the pair $(P_{d_i},S_{d_i})$ is observable and the  matrix $S_{d_i} \in \mathbb{R}^{n_{d_i} \times n_{d_i}}$ has the same properties as $S_r$ in \eqref{sigmodssg}. 
The global reference model \eqref{smodelg} and the local disturbance models \eqref{smodel} are merged into the \emph{signal model}
\begin{subequations}\label{sigmod}
\begin{align}
 \dot{w}(t) &= Sw(t), &&\hspace{-0.2cm} t > 0, \quad w(0) = w_{0} \in \mathbb{R}^{n_w} \label{sigmodssm}\\
 r(t) &= p\t w(t), &&\hspace{-0.2cm}  t \geq 0\\
 d_i(t) &= P_{i}w(t), &&\hspace{-0.2cm}  t \geq 0\label{sigmoy_routm}
\end{align}
\end{subequations}
so that all signal forms described in \eqref{smodelg} and \eqref{smodel} can be generated by \eqref{sigmod}. This directly determines the vector $p \in \mathbb{R}^{n_w}$ and the matrix $P_i \in \mathbb{R}^{m_i \times n_w}$, $i = 1,\ldots,N$. In the sequel, it is assumed that only $S$ in \eqref{sigmodssm} is known for the controller design. Note that $S$ inherits the properties from $S_r$ and $S_{d_i}$, i.\:e., $S$ is diagonalizable with eigenvalues on the imaginary axis.

The agents consist of two groups. The first group is composed of the agents $i$, $i = 1,\ldots,n$, $n \geq 1$, which have access to the reference input $r$ and are therefore called the \emph{informed agents}. In contrast, the information about the reference input can only be broadcast to the remaining agents $i$, $i = n+1,\ldots,N$, through a communication network. More precisely, these agents have only access to the reference information of their neighbours due to the communication constraints. Hence, they are the so-called \emph{uninformed agents}. As a consequence, a cooperative regulator is required, in order to achieve output regulation.

In this paper, the \emph{robust cooperative output regulation problem} is solved by utilizing
the \emph{cooperative state feedback regulator}
\begin{subequations}\label{imreg}
\begin{align}
 \!\!\!\dot{v}_i(t) &= Sv_i(t)\! +\! b_y\big(\sum_{j=1}^Na_{ij}(y_i(t)\! - \!y_j(t)) \!+ \!a_{i0}(y_i(t)\! -\! r(t))\big)\label{regODE}\\
  \!\!\!      u_i(t) &= u^l_i(t) + u^c_i(t) = \mathcal{K}_i[v_i(t),x(t)], \quad t \geq 0\label{sfeed}
\end{align}
\end{subequations}
for $i = 1,\ldots,N$ with \eqref{regODE} defined on $t > 0$, the IC $v_i(0) = v_{i,0} \in \mathbb{R}^{n_w}$ and $b_y \in \mathbb{R}^{n_w}$ such that the pair $(S,b_y)$ is controllable as well as $x = \col{x_1,\ldots,x_N} \in \mathbb{R}^{N}$ in \eqref{sfeed}. Furthermore, $\mathcal{K}_i$ is a formal \emph{feedback operator}, which is determined by the \emph{local state feedback}
\begin{equation}\label{lsf}
 u^l_i(t) = k_v\t v_i(t) - k_1x_i(1,t) -\int_0^1k_x(\zeta)x_i(\zeta,t)\d\zeta	 
\end{equation}
with the common feedback gains $k_v \in \mathbb{R}^{n_w}$, $k_1 \in \mathbb{R}$, $k_x(z) \in \mathbb{R}$ and the \emph{cooperative state feedback} 
\begin{equation}\label{gsfeed}
 u^c_i(t) = \int_0^1r_x(\zeta)\big(\sum_{j = 1}^Na_{ij}(x_i(\zeta,t) \! - \!x_j(\zeta,t)) \!+\! a_{i0}x_i(\zeta,t)\big)\d\zeta,
\end{equation}
where  $r_x(z) \in \mathbb{R}$ is the common feedback gain. 
\begin{rem}
In order to implement \eqref{gsfeed} with a low communication load, introduce $\xi_i(t) =  \int_0^1r_x(\zeta)x_i(\zeta,t)\d\zeta$ so that \eqref{gsfeed} can be rewritten as
\begin{equation}
 u^c_i(t) = \sum_{j = 1}^Na_{ij}(\xi_i(t) - \xi_j(t)) + a_{i0}\xi_i(t).
\end{equation}
Since all agents know the common feedback gain $r_x(z)$, they are able to transmit only the lumped quantities $\xi_i(t)$. \hfill $\triangleleft$
\end{rem}
In \eqref{regODE} and \eqref{gsfeed} the constants $a_{ij}$, $i,j = 1,\ldots,N$, are the elements of the adjacency matrix $A_{\mathcal{G}}$ corresponding to the digraph $\mathcal{G}$. By regarding \eqref{smodelg} as agent $0$, the constants $a_{i0}$, $i = 1,\ldots,N$, describe the communication between the reference model \eqref{smodelg} and the agents \eqref{drs}. More specifically, $a_{i0} > 0$ holds for the informed agents $i$, $i = 1,\ldots,n$, while $a_{i0} = 0$ is valid for the uninformed agents $i$, $i = n+1,\ldots,N$. In the following it is assumed that agent $0$ is the root of the digraph $\bar{\mathcal{G}}$ describing the communication network with node set $\bar{\mathcal{V}} = \{0,1,\ldots,N\}$ and edge set $\bar{\mathcal{E}}$. Then, by removing all edges of $\bar{\mathcal{E}}$, that are incident to the root, the subgraph $\mathcal{G}$ with node set $\mathcal{V} = \{1,\ldots,N\}$ and the edge set $\mathcal{E}$ is obtained. Note that the cooperative state feedback controller \eqref{gsfeed} has to respect the communication topology and thus leads to structural constraints, when designing the regulator \eqref{imreg}. The latter has to ensure stability of the networked controlled MAS and the \emph{reference tracking}
\begin{equation}\label{outregdualcond2}
\lim_{t \to \infty}e_{y_i}(t) = \lim_{t \to \infty}(y_i(t) - r(t)) = 0,
\end{equation}
$i = 1,\ldots,N$, for all ICs of the plant \eqref{drs}, of the signal model \eqref{sigmod} and of the controller \eqref{imreg}. Furthermore, the property \eqref{outregdualcond2} should be robust in the sense that it holds despite of all model uncertainties \eqref{modunc}, for which the nominal networked controller stabilizes the networked controlled MAS. This output regulation problem can also be seen as a \emph{leader-follower robust output synchronization problem} in the presence of disturbances, where the reference model \eqref{smodelg} is the leader and the agents \eqref{drs} are the followers.

An interesting specialization of the robust cooperative output regulation problem is  obtained by omitting the global reference model \eqref{smodelg}, i.\:e., by not specifying the reference input. This results in a \emph{leaderless robust output synchronization problem} subject to disturbances, where \eqref{outregdualcond2} is replaced by 
\begin{equation}\label{outregdualsync}
\lim_{t \to \infty}(y_i(t) - y_j(t)) = 0, \quad i,j = 1,\ldots,N.
\end{equation}
Here, the information exchange through the network is utilized so that the agents can negotiate a common synchronization trajectory for their outputs $y_i$. Different from the leader-follower output synchronization problem, the corresponding synchronization trajectory is not specified a priori, but depends on the network topology and on the ICs of the agents.

\section{Cooperative State Feedback Regulator Design}\label{sec:nomout}
The regulator \eqref{imreg} has a similar structure as the classical regulator resulting from the internal model principle (see \cite{Pa14,Deu19}). In particular, it also contains a copy of the signal model \eqref{sigmod}, but different from the usual approach it is driven by a diffusive coupling of neighbouring outputs, which takes the communication topology into account. This allows a solution of the considered output regulation problem by a suitable cooperation of the internal models for each agent. Therefore, \eqref{regODE} is called the \emph{cooperative internal model} in the sequel. 

In order to achieve robust cooperative output regulation, the state feedback controller \eqref{sfeed} must stabilize the nominal MAS. The latter results from setting the model uncertainties in \eqref{modunc} to zero. In addition, the cooperative internal model \eqref{regODE} is rewritten in the form
\begin{equation}\label{lapldintmod}
 \dot{v}_i(t) = Sv_i(t) + (e_i\t H \otimes b_y)y(t), \quad i = 1,\ldots,N,
\end{equation}
which results from a simple calculation by taking the \emph{leader-follower matrix}
\begin{equation}\label{Hdef}
 H = L_{\mathcal{G}} + L_0 \in \mathbb{R}^{N \times N}
\end{equation}
with $L_0 = \diag{a_{10},\ldots,a_{N0}}$ into account as well as denoting the $i$-th unit vector by $e_i \in \mathbb{R}^N$ and the Laplacian matrix w.r.t. the digraph $\mathcal{G}$ by $L_{\mathcal{G}}$. Furthermore, the reference input $r$ is not considered in \eqref{lapldintmod}, because this exogenous input does not influence the closed-loop stability (see \eqref{regODE}). With this, the definition of $v = \col{v_1,\ldots,v_N} \in \mathbb{R}^{Nn_w}$ directly leads to the \emph{aggregated cooperative internal model}
\begin{equation}\label{aggrdintmond}
 \dot{v}(t) = (I_N \otimes S)v(t) + (H \otimes b_y)y(t).
\end{equation}
Hence, the \emph{nominal networked controlled MAS} takes the form
\begin{subequations}\label{nomplant}
 \begin{align}
    \dot{v}(t) &= (I_N \otimes S)v(t) + (H \otimes b_y)\mathcal{C}[x(t)]\\
  \dot{x}(z,t) &= x''(z,t) + a(z)x(z,t)\\
  x'(0,t) &= q_0x(0,t)\\
  x'(1,t) &= q_1x(1,t) + \mathcal{K}[v(t),x(t)],\label{actbc}
 \end{align}	
\end{subequations}
in which the nominal output operator $\mathcal{C}$ results from \eqref{cdefindomain} by setting the model uncertainties \eqref{modunc} to zero. Furthermore, the definition $\mathcal{K} = \col{\mathcal{K}_1,\ldots,\mathcal{K}_N}$ was utilized in \eqref{actbc}. The stabilizing controller \eqref{sfeed} is obtained, by mapping \eqref{nomplant} into the stable \emph{ODE-PDE cascade}
\begin{subequations}\label{tsys}
\begin{align}
 \dot{e}_v(t) &= F_{e_v}e_v(t)\label{ODEtsys}\\
 \dot{\tilde{x}}(z,t) &= \tilde{x}''(z,t) - \mu_c\tilde{x}(z,t)\label{xtsyspde}\\
 \tilde{x}'(0,t) &= 0\\
 \tilde{x}'(1,t) &= (I_N \otimes k_v\t)e_v(t)\label{bctsys}
\end{align}	
\end{subequations}
with 
\begin{equation}\label{Fevdef}
 F_{e_v} = I_N \otimes S - H \otimes \tilde{q}(1)k\t_v.
\end{equation}
This requires to determine a \emph{local backstepping transformation}
\begin{equation}\label{btrafo}
 \tilde{x}(z,t) = x(z,t) - \int_0^zk(z,\zeta)x(\zeta,t)\d\zeta = \mathcal{T}_c[x(t)](z)
\end{equation}
with the common kernel $k(z,\zeta) \in \mathbb{R}$, which is applied to the individual agents \eqref{drs} without taking the communication into account. The latter is needed for the design of the \emph{cooperative decoupling transformation}
\begin{equation}
 e_{v_i}(t) = v_i(t)  - \int_0^1\!\!\tilde{q}(\zeta)\big(\sum_{j = 1}^Na_{ij}(\tilde{x}_i(\zeta,t) - \tilde{x}_j(\zeta,t)) + a_{i0}\tilde{x}_i(\zeta,t)\big)\d\zeta
\end{equation}
for $i = 1,\ldots,N$, which needs the state of several agents and describes the deviation $e_{v_i}$ of the ODE state $v_i(t)$ from the PDE states $x_i(z,t)$. Therein, the common vector $\tilde{q}(z) \in \mathbb{R}^{n_w}$ has to be determined. By introducing $e_v = \col{e_{v_1},\ldots,e_{v_N}} \in \mathbb{R}^{Nn_w}$ the aggregated decoupling transformation reads
\begin{equation}\label{decoupltrafo}
 e_{v}(t) = v(t) - \int_0^1(H \otimes \tilde{q}(\zeta))\tilde{x}(\zeta,t)\d\zeta,
\end{equation}
which follows from the same reasoning as for \eqref{aggrdintmond}. 

\subsection{Local Backstepping Transformation of the MAS}\label{sec:lback}
In the first step, the backstepping transformation \eqref{btrafo} is determined to map \eqref{nomplant} into the \emph{intermediate target system}
\begin{subequations}\label{tsyspre}
	\begin{align}
\!\!\!	\dot{v}(t) &= (I_N \otimes S)v(t) + (H \otimes b_y)\mathcal{C}\mathcal{T}_c^{-1}[\tilde{x}(t)]\label{vprefin}\\
\!\!\!	\dot{\tilde{x}}(z,t) &= \tilde{x}''(z,t) - \mu_c\tilde{x}(z,t)\label{pdes1}\\
\!\!\!	\tilde{x}'(0,t) &= 0\\
\!\!\!	\tilde{x}'(1,t) &= u(t) \! + \! \big(q_1 \!-\! k(1,1)\big)x(1,t) \! - \! \int_0^1\!\!\!\!k_z(1,\zeta)x(\zeta,t)\d\zeta.\label{BCtsysu} 
	\end{align}	
\end{subequations}
Therein, the agents are stabilized by choosing $\mu_c \in \mathbb{R}$ such that the PDE subsystem \eqref{xtsyspde}--\eqref{bctsys} is exponentially stable. Differentiating \eqref{btrafo} \wrt time and inserting \eqref{nomplant} and \eqref{tsyspre}, the same calculations as in \cite{Sm04} verify that $k(z,\zeta)$ has to solve the \emph{kernel equations}
\begin{subequations}\label{keq}
 \begin{align}
  k_{zz}(z,\zeta) - k_{\zeta\zeta}(z,\zeta) &= (\mu_c + a(\zeta))k(z,\zeta),\; 0 < \zeta < z < 1\\
  k_{\zeta}(z,0) &= q_0k(z,0)\\
  k(z,z) &= q_0 - \frac{1}{2}\int_0^z(\mu_c + a(\zeta))\d\zeta.	
 \end{align}
\end{subequations}
It is shown in \cite{Sm04} that \eqref{keq} has a unique $C^2$-solution. Furthermore, the inverse transformation exists and is given by
\begin{equation}\label{btrafoing}
 x(z,t) = \tilde{x}(z,t) + \int_0^zk_I(z,\zeta)\tilde{x}(\zeta,t)\d\zeta = \mathcal{T}^{-1}_c[\tilde{x}(t)](z).
\end{equation}
Therein, the kernel $k_I(z,\zeta) \in \mathbb{R}$ follows from similar kernel equations. In order to determine the operator $\mathcal{C}\mathcal{T}_c^{-1}$ in \eqref{vprefin}, insert \eqref{btrafoing} in \eqref{cdefindomain} for the nominal case and change the order of integration. This yields
\begin{equation}\label{ctcomp}
 \mathcal{C}\mathcal{T}_c^{-1}[\tilde{x}(t)] = \int_0^1\tilde{c}(\zeta)\tilde{x}(\zeta,t)\d\zeta + c_{b0}\tilde{x}(0,t) + c_{b1}\tilde{x}(1,t)
\end{equation}
with $\tilde{c}(\zeta) = c_{b1}k_I(1,\zeta) + c(\zeta) + \int_{\zeta}^1c(\bar{\zeta})k_I(\bar{\zeta},\zeta)\d\bar{\zeta}$ after straightforward computations.

\subsection{Decoupling of the Cooperative Internal Model}\label{sec:gdecoupl}
In the second step, the transformation \eqref{decoupltrafo} is utilized to map \eqref{tsyspre} into the final target system \eqref{tsys}. For this, the ODE subsystem \eqref{vprefin} is decoupled from the PDE subsystem \eqref{pdes1}--\eqref{BCtsysu}. In order to determine the state feedback \eqref{sfeed} in the original coordinates, the decoupling transformation \eqref{decoupltrafo} has to be represented in terms of $x$. Inserting \eqref{btrafo} in \eqref{decoupltrafo} and changing the order of integration results in
\begin{equation}\label{evorig}
 e_v(t) = v(t) - \int_0^1(H \otimes q(\zeta))x(\zeta,t)\d\zeta
\end{equation}
where $q(\zeta) = \tilde{q}(\zeta) - \int_{\zeta}^1\tilde{q}(\bar{\zeta})k(\bar{\zeta},\zeta)\d\bar{\zeta}$. With this, the state feedback \eqref{sfeed} is obtained from \eqref{BCtsysu} and by inserting \eqref{evorig} in \eqref{bctsys}. Particularly, the latter yields
\begin{equation}
 \tilde{x}'(1,t) = (I_N \otimes k_v\t)v(t) - \int_0^1(H \otimes k_v\t q(\zeta))x(\zeta,t)\d\zeta\label{bcnew}
\end{equation}
when taking the multiplication property of the Kronecker product into account (see Section \ref{sec:intro}). Consequently, the BC \eqref{BCtsysu} and \eqref{bcnew} yield the state feedback
\begin{multline}\label{ufeedfin}
 u(t) = (I_N \otimes k_v\t)v(t) - \big(q_1 - k(1,1)\big)x(1,t)\\
       + \int_0^1k_z(1,\zeta)x(\zeta,t)\d\zeta - \int_0^1 (H \otimes k_v\t q(\zeta))x(\zeta,t)\d\zeta. 
\end{multline}
Considering the aggregation of \eqref{lsf} and \eqref{gsfeed} and comparing the result with \eqref{ufeedfin} directly leads to the common feedback gains
\begin{equation}\label{nwfgains}
 k_1 = q_1 - k(1,1), k_x(\zeta) = -k_z(1,\zeta) \text{ and } r_x(\zeta) = -k_v\t q(\zeta).
\end{equation}

In order to  determine \eqref{decoupltrafo}, differentiate it \wrt time, insert \eqref{tsyspre} with \eqref{ctcomp} in the result and use the BC \eqref{bctsys}. After applying integrations by parts this results in
\begin{align}
\dot{e}_{v}(t) &= \dot{v}(t) 
                  - \int_0^1(H \otimes \tilde{q}(\zeta))\dot{\tilde{x}}(\zeta,t)\d\zeta\nonumber\\
               &= (I_N \otimes S - H \otimes \tilde{q}(1)k_v\t )e_{v}(t)\label{bvpder}\\    
& \quad + \big(H \otimes (b_yc_{b0} - \tilde{q}'(0))\big)\tilde{x}(0,t)\nonumber\\
& \quad + \big(H \otimes (\tilde{q}'(1) + b_yc_{b1})\big)\tilde{x}(1,t)\nonumber\\
& \quad + \int_0^1\!\!\!\big(H \otimes (b_y\tilde{c}(\zeta) - \tilde{q}''(\zeta) + S\tilde{q}(\zeta) + \mu_c\tilde{q}(\zeta))\big)\tilde{x}(\zeta,t)\d\zeta.\nonumber
\end{align}							
Hence, if $\tilde{q}(z)$ is the solution of the \emph{decoupling equations}
\begin{subequations}\label{bvp1}
	\begin{align}
	\tilde{q}''(z) - \mu_c\tilde{q}(z) - S\tilde{q}(z) &= b_y\tilde{c}(z), \quad z \in (0,1)\label{dode1}\\
	\tilde{q}'(0) &= b_y c_{b0}\label{dode2}\\
	\tilde{q}'(1) &= -b_yc_{b1},\label{dode3}
	\end{align}
\end{subequations}
then \eqref{decoupltrafo} and \eqref{ufeedfin} map the intermediate target system \eqref{tsyspre} into the final target system \eqref{tsys}. The next lemma asserts the solvability of \eqref{bvp1}. 
\begin{Lemma}[Solvability of the decoupling equations]\label{lem:decoupl}\hfill
 Denote the spectrum of the PDE subsystem \eqref{xtsyspde}--\eqref{bctsys} by $\sigma_c$. The \emph{decoupling equations} \eqref{bvp1} have a unique $C^2$-solution $\tilde{q}(z) \in \mathbb{R}^{n_w}$ if $\sigma_c \cap \sigma(S) = \emptyset$. 
\end{Lemma}
The proof of this lemma can be directly deduced from the corresponding result in \cite{Deu19}. Obviously, the condition of Lemma \ref{lem:decoupl} is fulfilled in the design, because $\sigma(S) \subset j\mathbb{R}$ and the PDE subsystem \eqref{xtsyspde}--\eqref{bctsys} is exponentially stable. Note that due to the constant coefficients in \eqref{bvp1}, the solution can be obtained explicitly. 

\subsection{Stability of the Networked Controlled MAS}
The resulting target system \eqref{tsys} shows that the constraint stabilization of the nominal networked controlled MAS \eqref{nomplant} originating from the restricted communication topology can be traced back to the constrained stabilization of the finite-dimensional ODE subsystem \eqref{ODEtsys}, \eqref{Fevdef}. This significantly facilitates the networked controller design for the infinite-dimensional MAS \eqref{drs}. In particular, systematic solutions exist for determining the common feedback gain $k_v\t$ in \eqref{Fevdef} to solve the \emph{simultaneous stabilization problem} resulting for \eqref{ODEtsys} (see, e.\:g., \cite[Ch. 5.4]{Is17} and \cite[Ch. 8.4]{Bu19}). For this it is required that the global reference inputs and the local disturbances can be transmitted from the agent input to its output. This leads to the nonblocking conditions for the numerator of the corresponding transfer behaviour in the next lemma.

\begin{Lemma}[Stabilization of the ODE subsystem]\label{lem:ctrl}\hfill The numerator $N(s) \in \mathbb{C}^{N\times N}$ of the transfer matrix $F(s) = N(s)D^{-1}(s)$ from $u$ to $y$ \wrt the nominal MAS \eqref{drs} is given by $N(s) = c_{b0}I_N + c_{b1} T\t \Psi\t(0,1,s) T 
+\int_0^1\tilde{c}(\zeta)T\t\Psi\t(0,\zeta,s)T\d\zeta$ with $T = [I_N \;\; 0]\t \in \mathbb{R}^{2N \times N}$ and
\begin{equation}\label{fmat}
\Psi(z,\zeta,s) = \e^{\begin{bmatrix}
	0              & I_N\\
	(s + \mu_c)I_N & 0
	\end{bmatrix}(z-\zeta)}
\end{equation}
(see \eqref{ctcomp}). Then, the pair $(S,\tilde{q}(1))$ is controllable iff the pair $(S,b_y)$ is controllable and $\det N(\lambda) \neq 0$, $\forall \lambda \in \sigma(S)$, holds. Assume that $(S,\tilde{q}(1))$ is controllable and let the digraph $\bar{\mathcal{G}}$ be connected with the node $0$ as its root, then there exists a common feedback gain $k_v\t$ such that $F_{e_v}$ in \eqref{Fevdef} is Hurwitz. A possible choice for this feedback  is 
\begin{equation}\label{obsgainrefnet}
 k_v\t = \tilde{q}\t(1)Q
\end{equation}
with $Q \in \mathbb{R}^{n_w \times n_w}$ the positive definite solution of the algebraic Riccati equation
\begin{equation}\label{are}
 S\t Q + QS - 2\nu Q\tilde{q}(1)\tilde{q}\t(1)Q + aI = 0,
\end{equation}
where $a > 0$ and $\nu$ such that $\operatorname{Re}\lambda \geq \nu > 0$, $\forall \lambda \in \sigma(H)$. 
\end{Lemma}
\begin{IEEEproof}
The calculation of the numerator $N(s)$ and the result for the controllability of the pair $(S,\tilde{q}(1))$ directly follows from the related result in \cite{Deu19}. If the digraph $\bar{\mathcal{G}}$ is connected with the node $0$ as its root, then $\operatorname{Re}\lambda > 0$, $\forall \lambda \in \sigma(H)$ (see \cite[Rem. 2]{Su13}) such that there exists a $\nu$ satisfying the condition of the lemma. With this,  and applying a similar reasoning as in \cite[Ch. 5.5]{Is17} it is easily verified that \eqref{obsgainrefnet} ensures a Hurwitz matrix $F_{e_v}$. In particular, consider the matrix $F_i = S - \lambda_i(H)\tilde{q}(1)k_v\t$, $i = 1,\ldots,N$, with $\lambda_i(H)$ being the eigenvalues of $H$ in \eqref{Hdef}. This matrix is a block diagonal element resulting from mapping $F_{e_v}$ to an upper triangular matrix, which is always possible (see \cite[Ch. 5.5]{Is17}). Consider $\xi^H(F_i^HQ + QF_i)\xi \leq \xi^H(S\t Q + QS - 2\nu Q\tilde{q}(1)\tilde{q}\t(1)Q)\xi = - a\|\xi\|^2$, $\xi \in \mathbb{C}^{n_w}$. Therein, \eqref{obsgainrefnet} is inserted in $F_i$ and the condition $\operatorname{Re}\lambda \geq \nu > 0$ as well as \eqref{are} are used (cf. \cite[Ch. 5.5]{Is17}). This and the fact that the algebraic Riccati equation \eqref{are} has a unique positive definite solution $Q$ for $(S,\tilde{q}(1))$ controllable implies that $F_i$ and thus $F_{e_v}$ are Hurwitz matrices. 
\end{IEEEproof}
\begin{rem}
It should be emphasized that the design only requires to solve the kernel equations \eqref{keq}, the decoupling equations \eqref{bvp1} and the Riccati equation \eqref{are}, that are independent from the number $N$ of agents. This verifies the \emph{scalability} of the proposed networked controller design. \hfill $\triangleleft$
\end{rem}
This completes the design of the networked controller \eqref{imreg}.   In the next theorem the stability of the resulting networked controlled nominal MAS is stated.
\begin{thm}[Nominal stability of the networked controlled MAS]\label{thm:clstab}\hfill
Let the feedback gains in \eqref{sfeed} be given by \eqref{nwfgains}. Assume that  $\mu_c > 0$ and that $F_{e_v}$  is a Hurwitz matrix such that $\alpha = \min(\alpha_{e_v},\mu_c) > 0$ where $\alpha_{e_v} = -\max_{\lambda \in \sigma(F_{e_v})}\operatorname{Re}\lambda$. Then, the abstract initial value problem (IVP) corresponding to the resulting nominal networked controlled MAS with the state $x_c(t) = \operatorname{col}(v(t),x(t))$ and $x(t) = \{x(z,t), z\in[0,1]\}$ is well-posed in the state space $X = \mathbb{C}^{Nn_w} \oplus (L_2(0,1))^N$ with the usual inner product. Furthermore, the system is exponentially stable in the norm 
$\|\cdot\|\; = (\|\cdot\|^2_{\mathbb{C}^{Nn_w}} +
\|\cdot\|_{L_2}^{2})^{1/2}$ where 
$\| h \|_{L_2}^2 = \int _ { 0 } ^ { 1 } \|h ( \zeta )\| _ { \mathbb { C } ^ { N } } ^ { 2 } \d \zeta$.
In particular, $\|x_c(t)\| \leq M\e^{-\alpha t}	
\|x_c(0)\|$, $t \geq 0$, holds for all $x_c(0) \in \mathbb{C}^{Nn_w} \oplus (H^2(0,1))^N \subset X$ satisfying the BCs of the closed-loop system and an $M \geq 1$. 
\end{thm}
\begin{IEEEproof}
Consider the transformation $\bar{x}(z,t) = \tilde{x}(z,t) - 0.5z^2(I_N \otimes k\t_v)e_v(t)$, which maps \eqref{tsys} into
\begin{subequations}\label{tsysp}
	\begin{align}
	\dot{e}_v(t) &= F_{e_v}e_v(t)\label{ODEtsysp}\\
	\dot{\bar{x}}(z,t) &= \bar{x}''(z,t) - \mu_c\bar{x}(z,t) + b\t(z)e_v(t)\label{xtsyspdep}\\
	\bar{x}'(0,t) &= 0\\
	\bar{x}'(1,t) &= 0\label{bctsysp}
	\end{align}	
\end{subequations}
for some $b(z) \in \mathbb{R}^{Nn_w}$. Define the operators $\mathcal{A}h = h'' - \mu_ch$, $h \in D(\mathcal{A}) = \{h \in (H^2(0,1))^{N} \,|\, h'(0) = h'(1) = 0\}$, $\mathcal{B}h = b\t h$, $h \in \mathbb{C}^{Nn_w}$ and $\tilde{\mathcal{A}}h = \col{F_{e_v}h_1,\mathcal{B}h_1 + \mathcal{A}h_2}$ with $h \in D(\tilde{\mathcal{A}}) = \mathbb{C}^{Nn_w} \oplus D(\mathcal{A})$ as well as introduce the states $\bar{x}(t) = \{\bar{x}(z,t), z \in [0,1]\}$ and $\xi(t) = \col{e_v(t),\bar{x}(t)}$. Then, \eqref{tsysp} can be represented by the abstract IVP $\dot{\xi}(t) = \tilde{\mathcal{A}}\xi(t)$, $t > 0$, $\xi(0) \in D(\tilde{\mathcal{A}}) \subset X$. Since $-\mathcal{A}$ is a \emph{Sturm-Liouville operator}, the operator $\mathcal{A}$ is the infinitesimal generator of a $C_0$-semigroup (see \cite{Del03,Deu19}). Then, the composite operator $\tilde{\mathcal{A}}$ is also an infinitesimal generator so that the abstract IVP in question is well-posed with the decay rate $\alpha$ in view of \cite[Lem. A-3.3]{Ha14}. In particular, the PDE subsystem \eqref{xtsyspdep}--\eqref{bctsysp} has decay rate $-\mu_c$, because its spectrum $\tilde{\sigma}$ satisfies $\max_{\lambda \in \tilde{\sigma}}\lambda \leq 0$ for $\mu_c = 0$ due to the Neumann BCs. By going through the chain of boundedly invertible transformations the stability result in the original coordinates can be verified.
\end{IEEEproof}

\section{Robust Cooperative Output Regulation}\label{sec:robout}
The networked controller \eqref{imreg} should be able to ensure reference tracking \eqref{outregdualcond2} in the presence of model uncertainties \eqref{modunc} and the modelled disturbances in \eqref{sigmod}. In order to investigate this property for the uncertain MAS \eqref{drs}, consider \eqref{regODE}
in the form
\begin{equation}\label{intmodr}
\dot{v}_i(t) = Sv_i(t) + b_ye_i\t( H y(t) - L_01_Nr(t)\big), \quad i = 1,\ldots,N.
\end{equation}
This result follows from a simple calculation using the properties of $L_{\mathcal{G}}$ and by defining $1_N = \col{1,\ldots,1} \in \mathbb{R}^N$. Observe that $H1_N = (L_{\mathcal{G}} + L_0)1_N = L_01_N$ holds (see \eqref{Hdef}), since  $L_{\mathcal{G}}1_N = 0$. With this and \eqref{intmodr}, the aggregated cooperative internal model \eqref{aggrdintmond} becomes
\begin{equation}\label{aggrdintmondr}
\dot{v}(t) = (I_N \otimes S)v(t) + (I_N \otimes b_y)\tilde{e}_y(t),
\end{equation}
in which
\begin{equation}\label{eytil}
 \tilde{e}_y(t) = He_y(t)
\end{equation}
and $e_y = y - 1_Nr$. The internal model \eqref{aggrdintmondr} has the same form as in \cite{Deu19} and thus fulfills the \emph{N-copy internal model principle}, which amounts to include $N$ copies of the signal model into the internal model. Hence, the corresponding results can be applied to verify robust cooperative output regulation. The next theorem makes this more precise.
\begin{thm}[Robust cooperative output regulation]\label{thm:robreg}
Assume that the model uncertainties \eqref{modunc} are such that the resulting networked controlled MAS is strongly asymptotically stable and the output operator $\bar{\mathcal{C}} = \col{\bar{\mathcal{C}}_1,\ldots,\bar{\mathcal{C}}_N}$ resulting from \eqref{cdefindomain} is relatively bounded. Furthermore, let the digraph $\bar{\mathcal{G}}$ be connected with the node $0$ as its root. Then, the networked controller \eqref{imreg} achieves robust cooperative output regulation, \ie $\lim_{t \to \infty}e_{y_i}(t) = 0$, $i = 1,\ldots,N$, independently of the disturbance input locations characterized by $g_{j,i}$, $j = 1,2,3,4$, $i = 1,\ldots,N$, as well as the generation of the disturbance and the reference signals by $p\t$, $P_i$ in \eqref{sigmod}.  
\end{thm}
\begin{IEEEproof}
The result in \cite{Deu19} ensures $\lim_{t \to \infty}\tilde{e}_y(t) = 0$ in the presence of the assumed model uncertainties \eqref{modunc} and exogenous signals generated by \eqref{sigmod}. If the digraph $\bar{\mathcal{G}}$ is connected with the node $0$ as its root, then $\det H \neq 0$ in view of \cite[Rem. 2]{Su13}. This and \eqref{eytil} imply $\lim_{t \to \infty}e_y(t) = 0$ verifying robust cooperative output regulation.
\end{IEEEproof}

\section{Leaderless Robust Output Synchronization}\label{sec:roboutsyncll}
In what follows, the results of the previous sections are applied to solve the \emph{leaderless robust output synchronization problem} introduced in Section \ref{sec:probform}. For this setup, the networked controller is given by \eqref{imreg} with $a_{i0} = 0$, $i = 1,\ldots,N$. Hence, the \emph{aggregated cooperative internal model} takes the form
\begin{equation}\label{aggrdintmonds}
\dot{v}(t) = (I_N \otimes S)v(t) + (L_{\mathcal{G}} \otimes b_y)y(t)
\end{equation}
in view of \eqref{Hdef} and \eqref{aggrdintmond}, because $L_0 = 0$. In this setup, the matrix $S_r$ in $S$ influences the form of the synchronization trajectory negotiated by the agents.  

\subsection{Backstepping Transformation into an ODE-PDE Cascade}
Similar to Sections \ref{sec:lback} and \ref{sec:gdecoupl}, the transformations
\begin{subequations}\label{synctrafcasc}
	\begin{align}
   \tilde{x}(z,t) &= x(z,t) - \int_0^zk(z,\zeta)x(\zeta,t)\d\zeta\label{btrafs}\\
         e_{v}(t) &= v(t) - \int_0^1(L_{\mathcal{G}} \otimes \tilde{q}(\zeta))\tilde{x}(\zeta,t)\d\zeta\label{decoupls}
	\end{align}
\end{subequations}
are utilized to map the nominal MAS \eqref{drs} into the ODE-PDE cascade \eqref{tsys} with
\begin{equation}\label{Fevs}
 F_{e_v} = I_N \otimes S - L_{\mathcal{G}} \otimes \tilde{q}(1)k\t_v.
\end{equation}
For this, the kernel $k(z,\zeta)$ in \eqref{btrafs} is the solution of \eqref{keq} and $\tilde{q}(\zeta)$ in \eqref{decoupls} has to result from solving \eqref{bvp1}. 

\subsection{Stabilization of the ODE Subsystem}
A somewhat different approach is needed to stabilize the ODE subsystem \eqref{ODEtsys} and \eqref{Fevs}, because $H$ in $F_{e_v}$ is replaced by the Laplacian matrix $L_{\mathcal{G}}$, which different from $H$ has an eigenvalue at the origin. As a consequence, the ODE subsystem has to be mapped into a cascade of two ODEs, in order to determine the feedback gain $k_v^T$. To this end, introduce the new  coordinates
\begin{equation}\label{syncerr}
  \col{e_{v_1}(t),
  \varepsilon_{v_2}(t),
  \ldots,
  \varepsilon_{v_N}(t)}
 = (\Theta \otimes I_{n_v})e_v(t)
\end{equation}
(see, e.\:g., \cite[Ch. 5.4]{Is17}) with $\varepsilon_{v_i} = e_{v_i} - e_{v_1}$, $i = 2,\ldots,N$, $\varepsilon_v = \col{\varepsilon_{v_2},\ldots,\varepsilon_{v_N}}$ and 
\begin{equation}\label{Thetadef}
\Theta = \begin{bmatrix}
1        & 0_{1 \times N-1}\\
-1_{N-1} & I_{N-1}
\end{bmatrix}.
\end{equation}
For the following, observe that
\begin{equation}\label{Ttraf}
\tilde{L}_{\mathcal{G}} = \Theta L_{\mathcal{G}}\Theta^{-1} = \begin{bmatrix}
0 & \tilde{l}\t_{12}\\
0_{N-1} & \tilde{L}_{22}
\end{bmatrix}
\end{equation}
holds, in which $0_{N-1} = \col{0,\ldots,0} \in \mathbb{R}^{N-1}$ and the inverse $\Theta^{-1}$ always exists. Then, applying the transformation \eqref{syncerr} to \eqref{tsys} and \eqref{Fevs} yields
\begin{subequations}\label{tsyss}
	\begin{align}
	\dot{\tilde{x}}(z,t) &= \tilde{x}''(z,t) - \mu_c\tilde{x}(z,t)\label{xtsyspdes}\\
      	 \tilde{x}'(0,t) &= 0\\
	     \tilde{x}'(1,t) &= 1_Nk_v\t e_{v_1}(t) + B\varepsilon_v(t)\label{bctsyss}\\
        \dot{e}_{v_1}(t) &= Se_{v_1}(t) - \big(\tilde{l}\t_{12} \otimes \tilde{q}(1)k\t_v\big)\varepsilon_v(t)\label{vlead}\\
    \dot{\varepsilon}_v(t) &= F_{\varepsilon_v}\varepsilon_v(t),\label{epssys}
   	\end{align}	
\end{subequations}
where
\begin{subequations}
\begin{align}
B &= \begin{bmatrix}
0\\
I_{N-1} \otimes k_v\t
\end{bmatrix}\\
F_{\varepsilon_v} &= I_{N-1} \otimes S - \tilde{L}_{22} \otimes \tilde{q}(1)k\t_v\label{hurmats}
\end{align}	 
\end{subequations}
after a simple calculation. In the following, $k_v\t$ is determined to solve the simultaneous stabilization problem arising for \eqref{epssys} (cf. \eqref{ODEtsys}). The next lemma is a direct consequence of the similarity between \eqref{Fevdef} and the matrix in \eqref{epssys} so that the result of Lemma \ref{lem:ctrl} becomes applicable.
\begin{Lemma}[Stabilization of the ODE subsystem]\label{lem:ctrls}\hfill Assume that $(S,\tilde{q}(1))$ is controllable (see Lemma \ref{lem:ctrl} for a condition) and let the digraph $\mathcal{G}$ associated with the Laplacian matrix $L_{\mathcal{G}}$ be connected, then there exists a common feedback gain $k_v\t$ such that the matrix $F_{\varepsilon_v}$ in \eqref{hurmats} is Hurwitz. A possible choice for this feedback is given by \eqref{obsgainrefnet} after solving \eqref{are} with $a > 0$ and $\nu$ such that
\begin{equation}\label{arebed}
 \operatorname{Re}\lambda \geq \nu > 0, \quad \forall \lambda \in \sigma(\tilde{L}_{22}).
\end{equation}
\end{Lemma}
\begin{IEEEproof}
It remains to verify the existence of $\nu$ such that \eqref{arebed} is satisfied. Since the digraph $\mathcal{G}$ is assumed to be connected, the Laplacian matrix $L_{\mathcal{G}}$ has only one eigenvalue at the origin and all other eigenvalues have a positive real part (see, e.\:g., \cite[Th. 5.1]{Is17}). Then, by \eqref{Ttraf} the matrix $\tilde{L}_{22}$ has only eigenvalues with positive real parts, which proves the lemma. 
\end{IEEEproof}
In order to investigate the stability properties of \eqref{tsyss}, introduce the transformations
\begin{subequations}\label{stabtrafo}
 \begin{align}
 \tilde{e}_{v_1}(t) &= e_{v_1}(t) - \Pi\varepsilon_v(t)\\
            \bar{x}(z,t) &= \tilde{x}(z,t) - \Sigma_1(z)\tilde{e}_{v_1}(t) - \Sigma_2(z)\varepsilon_v(t)\label{stabtrafo2}
 \end{align}
\end{subequations}
with $\Pi \in \mathbb{R}^{n_w \times (N-1)n_w}$, $\Sigma_1(z) \in \mathbb{R}^{N \times n_w}$ and $\Sigma_2(z) \in \mathbb{R}^{N \times (N-1)n_w}$. They map \eqref{tsyss} into
\begin{subequations}\label{tsyss2}
 \begin{align}
 \dot{\bar{x}}(z,t) &= \bar{x}''(z,t) - \mu_c\bar{x}(z,t)\label{xtsyspdes2}\\
 \bar{x}'(0,t) &= 0\\
 \bar{x}'(1,t) &= 0\label{bctsyss2}\\
 \dot{\tilde{e}}_{v_1}(t) &= S\tilde{e}_{v_1}(t)\label{vlead2}\\
 \dot{\varepsilon}_v(t) &= F_{\varepsilon_v}\varepsilon_v(t)\label{epssys2}
 \end{align}	
\end{subequations}
if  $\Pi$, $\Sigma_1(z)$ and $\Sigma_2(z)$ are the solution of
\begin{subequations}\label{kstabeq}
 \begin{align}
 \Pi F_{\varepsilon_v} - S\Pi &= -\big(\tilde{l}\t_{12} \otimes \tilde{q}(1)k\t_v\big)\\
 \Sigma_1''(z) - \mu_c\Sigma_1(z) - \Sigma_1(z)S &= 0, \quad z \in (0,1)\\
 \Sigma_1'(0) &= 0\\
 \Sigma_1'(1) &= 1_Nk_v\t\\
 \Sigma_2''(z) - \mu_c\Sigma_2(z) - \Sigma_2(z)F_{\varepsilon_v} &= 0, \quad z \in (0,1)\\
 \Sigma_2'(0) &= 0\\
 \Sigma_2'(1) &= 1_Nk_{v}\t\Pi + B.
 \end{align}
\end{subequations}
The next lemma clarifies the solvability of \eqref{kstabeq}.
\begin{Lemma}\label{lem:stabtraf} Assume that $\sigma(F_{\varepsilon_v}) \cap \sigma(S) = \emptyset$, $\sigma_c \cap \sigma(S) = \emptyset$ and $\sigma_c \cap \sigma(F_{\varepsilon_v}) = \emptyset$, in which $\sigma_c$ is the eigenvalue spectrum of the PDE subsystem \eqref{xtsyspdes}--\eqref{bctsyss}. Then, there exist unique matrices $\Pi \in \mathbb{R}^{n_w \times (N-1)n_w}$, $\Sigma_1(z) \in \mathbb{R}^{N \times n_w}$ and $\Sigma_2(z) \in \mathbb{R}^{N \times (N-1)n_w}$ solving \eqref{kstabeq}, where the elements of $\Sigma_i$, $i = 1,2$, are $C^2$-functions.
\end{Lemma}
The proof of this lemma directly follows from the related results in \cite{Deu19}. Note that the conditions of Lemma \ref{lem:stabtraf} can always be ensured by a suitable choice of $\mu_c$ (cf. Theorem \ref{thm:clstab}) and the design of the gain $k_v\t$ (cf. Lemma \ref{lem:ctrls}). The next theorem clarifies the stability properties of \eqref{tsyss}.
\begin{thm}[Nominal stability for output synchronization]\label{thm:outreg}
Assume that the conditions of Lemma \ref{lem:stabtraf} are fulfilled and that the PDE subsystem \eqref{xtsyspdes2}--\eqref{bctsyss2} is exponentially stable as well as let $F_{\varepsilon_v}$ be a Hurwitz matrix. Then, the solution of \eqref{tsyss} is bounded in the $L_2$-norm.
\end{thm} 
\begin{IEEEproof}
The assumptions of Theorem \ref{thm:outreg} imply $\lim_{t \to \infty}\|\bar{x}(t)\|_{L_2} = 0$ and $\lim_{t \to \infty}\varepsilon_v(t) = 0$. Then, in view of $\bar{x}(z,t) = \tilde{x}(z,t) - \Sigma_1(z)\tilde{e}_{v_1}(t) - \Sigma_2(z)\varepsilon_v(t)$ (see \eqref{stabtrafo2}), the solution $\tilde{x}(z,t)$ converges to 
\begin{equation}\label{xsync}
 \tilde{x}_{\infty}(z,t) = \Sigma_1(z)\tilde{e}_{v_1}(t).
\end{equation}
Hence, the solution of \eqref{tsyss} is bounded in the $L_2$-norm, because \eqref{vlead2} has a bounded solution by assumption (see Section \ref{sec:probform}).
\end{IEEEproof}
By inserting \eqref{xsync} in \eqref{ctcomp}, the steady state response
\begin{equation}\label{yinf}
 y_{\infty}(t) = \mathcal{C}\mathcal{T}^{-1}[\tilde{x}_{\infty}(t)] = \mathcal{C}\mathcal{T}^{-1}[\Sigma_1]\tilde{e}_{v_1}(t)
\end{equation}
is obtained. Therein, the evaluation of $\tilde{x}_{\infty}(z,t)$ is well-defined due to the smoothness of $\Sigma_1$ (see Lemma \ref{lem:stabtraf}). In what follows, it is verified that the elements of $y_{\infty}$ coincide in the presence of disturbances and model uncertainty, i.\:e., robust output synchronization is achieved. This highlights the fact that $S$ determines the form of the synchronization trajectory (see \eqref{vlead2}) and \eqref{yinf}.

\subsection{Robust Output Synchronization}
In order to verify robust output synchronization, the $N$-copy internal model principle has to hold for the cooperative internal model \eqref{aggrdintmonds}. From \eqref{aggrdintmondr} and \eqref{eytil} the result
\begin{equation}\label{aggrdintmonds2}
\dot{v}(t) = (I_N \otimes S)v(t) + (I_N \otimes b_y)L_{\mathcal{G}}y(t)
\end{equation}
is readily deduced in view of \eqref{Hdef},  $L_0 = 0$ and $a_{i0} = 0$, $i = 1,\ldots,N$. With $L_{\mathcal{G}} = \Theta^{-1}\tilde{L}_{\mathcal{G}}\Theta$ (see \eqref{Ttraf}) the result
\begin{equation}
 L_{\mathcal{G}}y(t) = \Theta^{-1}\tilde{L}_{\mathcal{G}}\begin{bmatrix}
 y_1(t)\\
 e_y(t)
 \end{bmatrix}
 = \tilde{H}e_y(t)
\end{equation}
follows, where $e_y = \col{y_2-y_1,\ldots,y_N-y_1}$ and
\begin{equation}\label{Htildef}
 \tilde{H} = \Theta^{-1}
 \begin{bmatrix}
  \tilde{l}\t_{12}\\
  \tilde{L}_{22}
 \end{bmatrix} \in \mathbb{R}^{N \times N-1}.
\end{equation}
In the proof of Lemma \ref{lem:ctrls} it is verified that a connected digraph $\mathcal{G}$ implies $\det \tilde{L}_{22} \neq 0$. Hence, $\operatorname{rank}\tilde{H} = N-1$ follows from \eqref{Htildef}. This verifies the $N$-copy internal model principle (see \cite{Deu19}), because $\tilde{H}e_y(t) = 0$ holds only for $e_y(t) = 0$. Consequently, $\lim_{t \to \infty}e_y(t) = 0$ and thus \eqref{outregdualsync} is ensured for non destabilizing model uncertainty and for the modeled disturbances. This is the result of the next theorem.
\begin{thm}[Leaderless robust output synchronization]\label{thm:robreg}\hfill
Assume that the model uncertainties \eqref{modunc} are such that the PDE subsystem (cf. \eqref{xtsyspdes}--\eqref{bctsyss} in the nominal case) resulting from applying the transformations \eqref{synctrafcasc} to the uncertain MAS \eqref{drs} and the cooperative internal model \eqref{aggrdintmonds}   is strongly asymptotically stable and the output operator $\bar{\mathcal{C}} = \col{\bar{\mathcal{C}}_1,\ldots,\bar{\mathcal{C}}_N}$ resulting from \eqref{cdefindomain} is relatively bounded. Furthermore, let the digraph $\mathcal{G}$ associated with the Laplacian matrix $L_{\mathcal{G}}$ be connected. Then, the networked controller \eqref{aggrdintmonds}, \eqref{sfeed} achieves robust output synchronization, \ie $\lim_{t \to \infty}(y_i(t) - y_j(t)) = 0$, $i,j = 1,\ldots,N$, independently of the disturbance input locations characterized by $g_{j,i}$, $j = 1,2,3,4$, $i = 1,\ldots,N$, and the generation of the disturbance by $P_i$ in \eqref{sigmod}.  
\end{thm} 
\begin{IEEEproof}
The result \eqref{vlead}--\eqref{epssys} shows that the transformed cooperative internal model can be split into stable subsystem \eqref{epssys} and the subsystem \eqref{vlead} with a bounded solution (see Section \ref{sec:probform}). With this, the result of \cite{Deu19} is directly applicable, since the dynamics of the latter subsystem can be merged with the dynamics of the disturbance model \eqref{smodel}. Consequently, robust output synchronization can be verified with the same reasoning as in \cite{Deu19}.
\end{IEEEproof}

\section{Example}
In order to demonstrate the results of the paper, consider $N = 4$ parabolic agents with the nominal parameters $a(z) = z+1$, $q_0 =3$ and $q_1 = 0$. The output to be controlled is determined by $c_0(z) = -z$, $c_{b0} = c_{b1} = 1$ and no pointwise in-domain measurement. Furthermore, the agents are affected by the local constant disturbances $d_i(t) = d_{i0}\in \mathbb{R}$, $i = 1,\ldots,4$, with the same disturbance model $S_{d_i} = S_d = 0$ and $p_{d_i} = 1$. They act at the agents according to the disturbance input locations $g_{1,1}(z) = 2z$, $g_{1,2}(z) = 3z+1$, $g_{1,3}(z) = z-1$, $g_{1,4}(z) = 2z$, $g_{2,i} = g_{3,i} = 1$ and $g_{4,i} = 0$, $i = 1,\ldots,4$. The reference model coinciding with the leader is denoted as agent $0$ and generates a sinusoidal reference output $r(t) = A_r\sin(\pi t + \varphi_r)$, $A_r,\varphi_r \in \mathbb{R}$. This leads in \eqref{smodelg} to
\begin{equation}
 S_r = \begin{bmatrix}
          0 & \pi\\
       -\pi & 0
 \end{bmatrix}, \quad p_r\t = \begin{bmatrix}
 1 & 0
 \end{bmatrix}
\end{equation}
so that $S = \operatorname{bdiag}(S_r,S_d)$ holds in \eqref{regODE}.
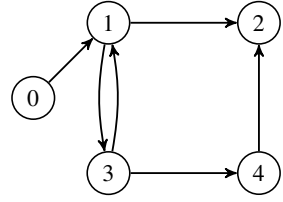
\begin{figure}
	\begin{center}
		\begin{minipage}{0.5\linewidth}
			\centering
			\begin{tikzpicture}[->,>=stealth',shorten >=0pt,auto,node distance=2.8cm,
			semithick,scale=2]
			\node (0) at (0,0.5) [circle,draw] {0};
			\node (1) at (0.5,1) [circle,draw] {1};
			\node (2) at (1.5,1) [circle,draw] {2};
			\node (3) at (0.5,0) [circle,draw] {3};
			\node (4) at (1.5,0) [circle,draw] {4};
			\draw[line width=0.7pt] (0) to (1);
			\draw[line width=0.7pt] (1) to[bend right=10] (3);
			\draw[line width=0.7pt] (3) to[bend right=10] (1);
			\draw[line width=0.7pt] (1) to (2);
			\draw[line width=0.7pt] (3) to (4);
			\draw[line width=0.7pt] (4) to (2);
			\end{tikzpicture}
		\end{minipage}
	\end{center}
\caption{Communication graph $\bar{\mathcal{G}}$ with agent $0$ (reference model) as leader.}\label{fig:graph}
\end{figure}%
The agents are able to transfer information through a communication network described by the digraph $\bar{\mathcal{G}}$ in Figure \ref{fig:graph} with the Laplacian matrix
\begin{equation}
 L_{\mathcal{G}} = \begin{bmatrix} 1 & 0 & -1 & 0\\ 
 -1 & 2 & 0 & -1\\
 -1 & 0 & 1 & 0\\  
 0  & 0 & -1 & 1\end{bmatrix}
\end{equation}
for the subgraph $\mathcal{G}$ with node set $\mathcal{V} = \{1,2,3,4\}$ and $L_0 = \operatorname{diag}(1,0,0,0)$. It is not difficult to verify that $\bar{\mathcal{G}}$ and $\mathcal{G}$ are connected, in which the agent $0$ is the root of the former graph. 

For the design of the networked controller achieving robust cooperative output regulation, the vector $b_y = [1 \; 1\; 1]\t$ ensuring $(S,b_y)$ controllable is chosen. Subsequently, the kernel equations \eqref{keq} are solved for $\mu_c = -5$ with the method of successive approximations (see \cite{Sm04}). After solving the decoupling equations \eqref{bvp1} and verifying $\det N(\lambda) \neq 0$, $\lambda \in \{0,\pm\text{j}\pi\}$, a solution of the algebraic Riccati equation \eqref{are} is obtained for $\nu = 0.382$ and $a = 150$. 

\begin{figure}[t!]
	\centering
\centering
    \input{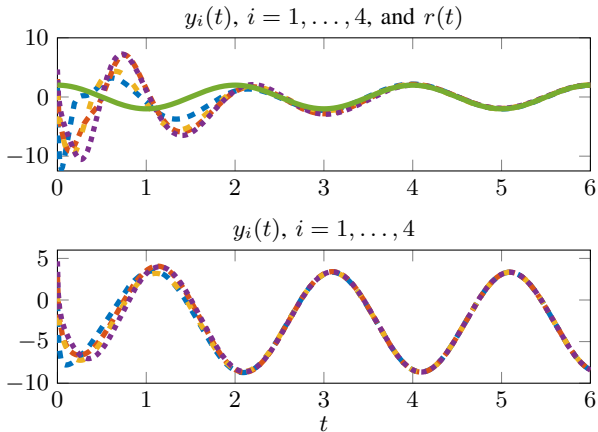}\vspace*{-1ex}
	\caption{Tracking behaviour for the agent $1$ (\ref{Agent1}),  agent $2$ (\ref*{Agent2}),  agent $3$ (\ref*{Agent3}) and agent $4$ (\ref*{Agent4}) in the presence of model uncertainty and local constant disturbances $d(t) = [3 \;\; -3 \;\; 1 \;\;1]\t$. Upper plot: tracking behaviour w.r.t. the leader output $r(t) = 2\cos(\pi t)$ (\ref*{Reference}); lower plot: synchronization behaviour without leader.}\label{Fig:sim}
\end{figure}

In order to verify robust cooperative output regulation, the model uncertainties $\Delta \lambda_1 = 0.2$, $\Delta \lambda_2 = -0.2$, $\Delta \lambda_3 = -0.1$, $\Delta \lambda_4 = 0.1$, $\Delta a_1(z) = 0.2a(z)$, $\Delta a_2(z) = -0.2a(z)$, $\Delta a_3(z) = 0.1a(z)$, $\Delta a_4(z) = 0.1a(z)$, $\Delta c_{0,i}(z) = 0$, $\Delta c_{b0,4} = -0.05$, $\Delta c_{b1,4} = 0.1$, and all other $\Delta c_{bk,i}$ vanishing are assumed. The resulting networked MAS is simulated for the ICs $w_r(0) = [2 \;\; 0]\t$, $x(z,0) = [1 \;\; 2 \;\; 0.5 \;\; 3]\t$, $v_1(0) = [1 \;\; 3.5 \;\; 0.5]\t$, $v_2(0) = [0.1 \;\; 2 \;\; 0.8]\t$, $v_3(0) = [1.7 \;\; 0.8 \;\; 0.3]\t$ and $v_4(0) = [0.5 \;\; 0.7 \;\; 0.9]\t$. The upper plot in Figure \ref{Fig:sim} shows the simulation results and verifies robust cooperative output regulation, i.\:e., the tracking of the leader output in the presence of different local disturbances $d(t) = [3 \;\; -3 \;\; 1 \;\;1]\t$ and model uncertainty. 

In order to investigate leaderless robust output synchronization, the networked controller \eqref{aggrdintmonds} and \eqref{sfeed} is applied to the same uncertain MAS. For this, also the same design parameters are utilized but $\nu = 1$ is chosen, in order to satisfy \eqref{arebed}. The resulting synchronization behaviour is depicted in the lower plot of Figure \ref{Fig:sim}. Obviously, the synchronization trajectory is particularly determined by the ICs of the agents, which can be seen from its amplitude matching the amplitudes of the agents ICs. Note that in the upper plot all agents outputs synchronize with the reference trajectory specified by the leader despite of disturbances and model uncertainty. This is different for the leaderless synchronization in the lower plot, where the synchronization trajectory depends also on the model uncertainty and the disturbances.

\section{Concluding remarks}
The results of the paper are directly extendible to obtain cooperative output feedback regulators by designing local backstepping observers for agents with boundary measurements. From the perspective of system classes further extensions concern parabolic agents with both temporal and spatially varying coefficients, coupled PDEs and hyperbolic agents. As far as the communication topology is considered, the inclusion of time-varying and random digraphs are of interest. 

%



\bibliographystyle{IEEEtranS}
\bibliography{IEEEabrv,mybib}

\end{document}